\documentclass[final,twoside,11pt]{entics}
\usepackage{enticsmacro}

\usepackage{stmaryrd}
\usepackage{mathtools}
\usepackage{amsmath}
\usepackage[capitalize]{cleveref}

\usepackage{quiver}

\newenvironment{eqalign}{\begin{equation}\begin{aligned}}{\end{aligned}\end{equation}}
\newenvironment{eqalign*}{\begin{equation*}\begin{aligned}}{\end{aligned}\end{equation*}}

\newenvironment{diagram*}{\begin{equation*}\begin{tikzcd}}{\end{tikzcd}\end{equation*}}

\tikzset{
  relation/.style={
    draw=none,
    every to/.append style={
      edge node={node [sloped, allow upside down, auto=false]{$#1$}}}
  }
}

\tikzcdset{
  mono/.code={
    \pgfsetarrows{tikzcd to reversed-tikzcd to}
  }
}
\tikzcdset{
  epi/.code={
    \pgfsetarrowsend{tikzcd double to}
  }
}
\tikzcdset{
  into/.code={
    \pgfsetarrows{tikzcd right hook-tikzcd to}
  }
}
\tikzcdset{
  twocell/.style={Rightarrow, shorten >= 3ex, shorten <= 3ex}
}

\tikzcdset{
  row sep/normal={
    6ex
  },
  column sep/normal={
    8ex
  }
}

\newcommand{\Overset}[2]{%
  \mathop{#2}\limits^{\vbox to -.1ex{%
  \kern -1.8ex\hbox{$#1$}\vss}}%
}
\newcommand{\Underset}[2]{%
  \mathop{#2}\limits_{\vbox to .1ex{%
  \kern -.6ex\hbox{$#1$}\vss}}%
}

\newcommand{\comp}{\fatsemi}

\mathchardef\dash="2D

\newcommand{\longto}{\longrightarrow}

\newcommand{\narrow}[2]{\overset{#1}{#2}}
\newcommand{\nto}[1]{\narrow{#1}{\to}}

\newcommand{\nlongto}[1]{\xrightarrow{#1}}

\newcommand{\opticto}{\leftrightarrows}
\newcommand{\chartto}{\rightrightarrows}

\newcommand{\de}{\mathrm{d}}

\newcommand{\biglens}[2]{
	 \begin{pmatrix}{\vphantom{f_f^f}#1} \\ {\vphantom{f_f^f}#2} \end{pmatrix}
}
\newcommand{\littlelens}[2]{
	 \begin{psmallmatrix}{\vphantom{f}#1} \\ {\vphantom{f}#2} \end{psmallmatrix}
}
\newcommand{\lens}[2]{
  \relax\if@display
	 \biglens{#1}{#2}
  \else
	 \littlelens{#1}{#2}
  \fi
}

\usepackage{xstring}
\newcommand{\cat}[1]{
  \relax
  \StrLen{#1}[\catarglen]
  \ifnum\catarglen=1
    \mathcal{#1}
  \else
    \mathbf{#1}
  \fi
}

\newcommand{\cod}{\mathrm{cod}}
\newcommand{\dom}{\mathrm{dom}}

\newcommand{\Sub}{\mathrm{Sub}}

\newcommand{\view}{\mathrm{view}}

\newcommand{\id}{\mathrm{id}}

\newcommand{\iso}[1][]{\overset{#1}{\cong}}
\newcommand{\equi}{\simeq}

\newcommand{\Cat}{\cat{Cat}}

\newcommand{\Set}{\cat{Set}}

\newcommand{\op}{\mathsf{op}}

\newcommand{\Para}{\cat{Para}}

\newcommand{\Lens}{\cat{Lens}}

\newcommand{\Fib}{\twocat{F}\cat{ib}}

\newcommand{\fflat}{{\flat\kern-1.4pt\flat}}

\renewcommand{\Fib}{\cat{Fib}}

\newcommand{\tower}[1]{\vec{#1}}
\newcommand{\Simple}[1]{\cat{S}(#1)}
\newcommand{\SubFib}[1]{{#1^\subseteq}}
\newcommand{\sub}{\mathrm{sub}}
\newcommand{\simple}{S}%{\mathrm{simp}}

\newcommand{\Fam}{\cat{Fam}}

\newcommand{\sDial}{\mathbf{Dial}_\times}
\newcommand{\Dial}{\mathbf{Dial}}
\newcommand{\DepDial}{\mathbf{Dial}}

\newcommand{\Sum}{\Fam}%{\mathbf{Sum}}
\newcommand{\Prod}{\mathbf{Cofam}}%{\mathbf{Prod}}
\newcommand{\sSum}{\mathbf{Sum}}%_\times}
\newcommand{\sProd}{\mathbf{Prod}}%_\times}

\newcommand{\car}{\mathrm{cart}}
\newcommand{\opcar}{\mathrm{opcart}}
\newcommand{\ver}{\mathrm{vert}}

\newcommand{\mono}{\mathsf{mono}}

\newcommand{\bigdial}[3]{
	 \begin{pmatrix}{\vphantom{f_f^f}#1} \\ {\vphantom{f_f^f}#2} \\ {\vphantom{f_f^f}#3} \end{pmatrix}
}
\newcommand{\littledial}[3]{
	 \begin{psmallmatrix}{\vphantom{f}#1} \\ {\vphantom{f}#2} \\ {\vphantom{f}#3} \end{psmallmatrix}
}
\newcommand{\dial}[3]{
  \relax\if@display
	 \bigdial{#1}{#2}{#3}
  \else
	 \littledial{#1}{#2}{#3}
  \fi
}

\DeclareFontFamily{U}{musix}{}%
\DeclareFontShape{U}{musix}{m}{n}{%
  <-12>   musix11
  <12-15> musix13
  <15-18> musix16
  <18-23> musix20
  <23->   musix29
}{}%
\newcommand*\musix{\usefont{U}{musix}{m}{n}\selectfont}
\DeclareTextFontCommand{\textmusix}{\musix}

\newcommand{\doublesharp}{{\raisebox{.6ex}{\textmusix{5}}}}

\newcommand{\dual}[2][]{{{#2}^{\vee_{#1}}}}

\renewcommand{\Fib}{\cat{Fib}}

\def\conf{MFPS 2024}
\def\myconf{mfps40} 	%%Fill in the acronym for your conference (with year)
\volume{4}			%Fill in the ENTICS volume number here
\def\volu{4}			% and here
\def\papno{7}			%Fill in your paper number here
\def\mydoi{14638}%
\def\lastname{Capucci, Gavranovi\'{c}, Malik, Rios \&\ Weinberger} %Lastnames appear in the running header 
\def\runtitle{On a fibrational construction\ldots} %Short title appears in the running header on even pages. 
\def\copynames{M. Capucci, B. Gravranovi\'{c}, A. Malik, }
\def\copyname{\ \ F. Rios, J. Weinberger\ }    %% Fill in the first initial and last name of the authors

\begin{document}
\begin{frontmatter}
\title{On a Fibrational Construction for\\[1ex]
Optics, Lenses, and Dialectica Categories}
\author{Matteo Capucci\thanksref{MCemail}}
\address{Mathematically Structured Programming, University of Strathclyde, Glasgow, UK}
\author{Bruno Gavranovi\'{c}\thanksref{BGemail}}
\address{Symbolica AI}
\author{Abdullah Malik\thanksref{AMemail}}
\address{Department of Mathematics, Florida State University, Tallahassee, FL, USA}
\author{Francisco Rios\thanksref{FRemail}}
\address{Quantum Computational Science Group, Oak Ridge National Laboratory, Oak Ridge, TN, USA}
\author{Jonathan Weinberger\thanksref{JWemail}}
\address{Fowler School of Engineering \& Center of Excellence in Computation, Algebra, and Topology (CECAT),\\ Chapman University, Orange, CA, USA}

\thanks[MCemail]{\href{mailto:matteo.capucci@strath.ac.uk}{matteo.capucci@strath.ac.uk}}
\thanks[BGemail]{\href{mailto:bruno@brunogavranovic.com}{bruno@brunogavranovic.com}}
\thanks[AMemail]{\href{mailto:amalik3@fsu.edu}{amalik3@fsu.edu}}
\thanks[FRemail]{\href{mailto:riosfr@ornl.gov}{riosfr@ornl.gov}}
\thanks[JWemail]{\href{mailto:jweinberger@chapman.edu}{jweinberger@chapman.edu}}

	\begin{abstract}
		Categories of lenses/optics and Dialectica categories are both comprised of bidirectional morphisms of basically the same form.
		In this work, we show how they can be considered a special case of an overarching fibrational construction, generalizing Hofstra's construction of Dialectica fibrations and Spivak's construction of generalized lenses.
		This construction turns a tower of Grothendieck fibrations into another tower of fibrations by iteratively twisting each of the components, using the opposite fibration construction. \\

%		\textbf{This is the peer-reviewed version of an accepted contribution to the 40th Conference on Mathematical Foundations of Programming Semantics 2024 (MFPS 2024), \url{https://oxford24.github.io/}, to be published in the proceedings.} \\
	\end{abstract}
\end{frontmatter}

	\section{Introduction}
	At various points in the past 40 years, people have come up with notions of `bidirectional morphisms' meant to represent two-ways interactions of various sorts.
	In this work we will be concerned with three specific kinds of bidirectional transformations.

	The first come from~\cite{de_paiva_dialectica_1989}, where de~Paiva defines \emph{Dialectica categories} as models of linear logic. The bidirectional morphisms therein represent an interactive proof in which parties exchange proofs and counterexamples.
	The second are lenses~\cite{foster_combinators_2007}, which have been introduced in the database theory community as a structure to encode view/update operations on record-type data structures.
	The third are optics~\cite{pickering_profunctor_2017,clarke_profunctor_2024}, which have been introduced to generalize lenses to a wider range of data structures.

	Coming from an applied category theory background, the authors have been mainly interested in bidirectional transformations to represent the compositional structure of systems as diverse as machines and games~\cite{Capucci_2022,myers2022categorical,niu2023polynomial}.

	A fairly complete recollection of the intricate history of these discoveries can be found in~\cite{hedges_lenses_2018}. It's interesting to notice the remarkable variety of motivations that have brought people to use or even rediscover lens-like structures.
	Nonetheless, aspects like \emph{lawfulness}, which are really important for lenses and optics as pattern for manipulating data structures, turn out to not matter in other settings, suggesting a fairly simple structure underlies all these examples, without necessarily exhausting them.

	In~\cite{spivak_generalized_2019} Spivak suggested that behind all the uses of lenses there is a very elementary construction in fibered category theory, the \emph{dual fibration} construction.
	This is the fibered variant of the `opposite category' construction: it takes a fibered category $P:\cat E \to \cat X$ and yields a new fibration $\dual P: \dual{\cat E} \to \cat X$ with the property that the fibers of the latter are the opposite of the fibers\footnote{See Definition~\ref{def:fiber} below for a precise definition.} of the former:
	\begin{equation*}
		\de {\dual P}(X) \equi \de P(X)^\op, \qquad X \in \cat X.
	\end{equation*}
	We write $\lens{E}{X}$ for an object $E : \cat E$ such that $P(E) =X$. In fact, fibering a category $\cat E$ means also separating its morphisms $\lens{E}{X} \to \lens{E'}{X'}$ into a \emph{cartesian} part $f:X \to X'$, living in the base $\cat X$, and a \emph{vertical} part $f^\flat : E \to f^*E'$, living in the fiber $\de P(X)$.

	Notice the cartesian part can be defined independently of the vertical part, by simply giving a morphism in $\cat X$, whereas the vertical part can't.
	This induces a natural `arrow of time', given by the order $f$ and $f^\flat$ have to be defined in:
	\begin{equation*}
		\lens{f^\flat}{f} : \lens{E}{X} \chartto \lens{E'}{X'}, \quad
		f:X \to X'  :  \cat X,\ f^\flat : E \to f^*E'  :  \de P(X)
	\end{equation*}
	Morphisms in $\dual{\cat E}$ are still bipartite in the same fashion, but now the vertical part points in the opposite direction than the cartesian one.
	Therefore morphism in $\dual{\cat E}$ seem to describe two consecutive `rounds of interaction' between two parties, with the characteristic back and forth of lens-like morphisms:
	\begin{equation*}
		\lens{f^\sharp}{f} : \lens{E}{X} \opticto \lens{E'}{X'}, \quad
		f:X \to X'  :  \cat X,\ f^\sharp : f^*E' \to E  :  \de P(X)
	\end{equation*}

	This pattern captures many examples of bidirectional morphisms in literature (from wiring diagrams to morphisms of locally ringed spaces), but fails to capture other lens-like morphisms such as optics and Dialectica morphisms.

	In this paper we reconcile these three classes of examples into one by extending the construction of duals to \emph{towers of fibrations}, \emph{i.e.}~sequences of fibrations:
	\begin{equation*}
		\tower{P} = \cat E_n \nlongto{P_n} \cdots \longto \cat E_2 \nlongto{P_2} \cat E_1 \nlongto{P_1} \cat E_0
	\end{equation*}
	Towers like these display $n$ consecutive rounds of interactions between two parties. By iteratively dualize each fibration in the sequence, we can then produce a category $\dual[n]{\cat E_n}$, still fibered over $\cat E_0$ (and over suitable duals of the intermediate categories) whose morphisms are made of $n$ parts alternating in direction. We call these \textbf{dialenses}, as a \emph{portmanteu} of `Dialectica' and `lenses'.

	Indeed, we show that both optics and Dialectica morphisms are dialenses of height $2$, and obviously lenses are dialenses of height $1$.

	\paragraph{Outline of the paper}
	In Section~\ref{sec:fibs} we recall the necessary notions of fibred category theory, as well as introduce the main examples of fibrations we use later.
	In Section~\ref{sec:examples} we recall the definitions of (generalized) lenses, Dialectica morphisms, and optics.
	In Section~\ref{sec:dialenses} we introduce the iterated dual construction and show how they capture the examples of interest.
	Finally, in~Section~\ref{sec:hofstra} we compare our construction of Dialectica categories with that of Hofstra, who also used fibrational tools.

	\paragraph{Notation}
	We denote by $f\comp g$ the composition of morphisms $f$ and $g$ in diagrammatic order.
	We make generous use of generalized elements, written $x:X$, to talk about objects in categories.

	\section{Preliminary notions on fibrations}
	\label{sec:fibs}

	\begin{definition}(\textbf{Fiber of a functor})\label{def:fiber}
		Let $P: \cat E \to \cat X$ be a functor. For an object $X  :  \cat X$, we define the \emph{fiber} of $P$ at $X$ as the subcategory $\de P(X)$ of $\cat E$ consisting of those morphisms $\varphi$ in $\cat E$ such that $P(\varphi) = \id_X$, which we call \emph{vertical}.
	\end{definition}

	The notation $\de P(X)$ to denote the fiber of $P$ at $X$ is due to B\'enabou~\cite{benabou_distributors_2000}.
	Other popular choices are $P^{-1}(X)$ and $\cat E_X$.

	\begin{definition}(\textbf{Cartesian arrow})\label{def:cart-arr}
		Let $P: \cat E \to \cat X$. A morphism $\varphi : E \to D$ in $\cat E$ is \emph{(strong) $P$-cartesian} if it satisfies the following universal property: let us write $f := P(\varphi)$ and $f : X \to Y$ in $\cat X$. Then we demand that for any $g : Z \to X$ in $\cat X$ and any $g' : U \to D$ in $\cat E$ with $P(g') = g \comp f$, there exist a unique morphism $h : U \to E$ in $\cat E$ with $P(h) = g$ and $h \comp \varphi = g'$:
		\begin{equation}
			\begin{tikzcd}[ampersand replacement=\&, sep=4ex]
				U \\
				Z \& E \& \& D \\
				\& X \& \& Y
				\arrow[dotted, maps to, from=2-2, to=3-2]
				\arrow["{{\exists ! h}}", dashed, from=1-1, to=2-2]
				\arrow[dotted, maps to, from=1-1, to=2-1]
				\arrow["{{\forall g}}"', curve={height=6pt}, from=2-1, to=3-2]
				\arrow["f", from=3-2, to=3-4]
				\arrow["\varphi", from=2-2, to=2-4]
				\arrow[dotted, maps to, from=2-4, to=3-4]
				\arrow["{{\forall g'}}"{description}, curve={height=-24pt}, from=1-1, to=2-4]
			\end{tikzcd}
		\end{equation}
	\end{definition}

	\begin{definition}(\textbf{Grothendieck fibration})
		A \emph{Grothendieck fibration} (or \emph{fibration} for short) is a functor $P:\cat E \to \cat X$ such that any arrow $f:X \to P(D)$ in $\cat X$ admits a (strong) cartesian lift:
		\begin{equation}
			% https://q.uiver.app/#q=WzAsNixbMSwxLCJYIl0sWzIsMSwiUChWKSJdLFsyLDAsIlYiXSxbMSwwLCJmXipWIl0sWzAsMSwiXFxjYXQgQiJdLFswLDAsIlxcY2F0IEUiXSxbMywyLCJmX1YiXSxbMCwxLCJmIl0sWzUsNCwiUCIsMl1d
			\begin{tikzcd}[ampersand replacement=\&]
				{\cat E} \& {f^*D} \& D \\
				{\cat X} \& X \& {P(D)}
				\arrow["{f_D}", from=1-2, to=1-3]
				\arrow["f", from=2-2, to=2-3]
				\arrow["P"', from=1-1, to=2-1]
			\end{tikzcd}
		\end{equation}
		For $f_D$ to be a cartesian lift, it has to be terminal among all morphisms into $D$ and projecting down to $f$, see Definition~\ref{def:cart-arr}. Moreover, we can show these lifts compose: cartesian lifts of $f$ and $g$ compose to give a cartesian lift for $f \comp g$. Furthermore, up to unique vertical isomorphism, the cartesian lift $f_D : f^*D \to D$ is determined uniquely.
		% This can be expressed by saying there exists a right adjoint $\ell$ in $\Cat/B$ to the map sending an object $V:\cat E$ to its identity morphism $P(V) \equalto P(V) : \cat X/P$:
		% \begin{equation}
		% 	% https://q.uiver.app/#q=WzAsMyxbMCwwLCJcXGNhdCBCL1AiXSxbMiwwLCJcXGNhdCBFIl0sWzEsMSwiXFxjYXQgQiJdLFsxLDAsIlxcaWQiLDIseyJjdXJ2ZSI6MX1dLFswLDEsIlxcZWxsIiwyLHsiY3VydmUiOjF9XSxbMCwyLCJcXHBpIiwyXSxbMSwyLCJQIl0sWzMsNCwiIiwyLHsibGV2ZWwiOjEsInN0eWxlIjp7Im5hbWUiOiJhZGp1bmN0aW9uIn19XV0=
		% 	\begin{tikzcd}[ampersand replacement=\&]
		% 		{\cat X/P} \&\& {\cat E} \\
		% 		\& {\cat X}
		% 		\arrow[""{name=0, anchor=center, inner sep=0}, "\id"', curve={height=6pt}, from=1-3, to=1-1]
		% 		\arrow[""{name=1, anchor=center, inner sep=0}, "\ell"', curve={height=6pt}, from=1-1, to=1-3]
		% 		\arrow["\pi"', from=1-1, to=2-2]
		% 		\arrow["P", from=1-3, to=2-2]
		% 		\arrow["\dashv"{anchor=center, rotate=-90}, draw=none, from=0, to=1]
		% 	\end{tikzcd}
		% \end{equation}

		$\cat E$ is called the \emph{total category} and $\cat X$ the \emph{base category}.
	\end{definition}

	We are going to assume that every fibration is \emph{cloven}, i.e.~comes equipped with a specific choice of cartesian lifts. A way to say this is that $\de P$ is an indexed category (\emph{i.e.}~a well-defined pseudofunctor $\cat X^\op \to \Cat$) and we are given an explicit isomorphism $\int \de P \iso P$.
	In this way we can move freely between a fibration and its indexing of fibers.

	\begin{proposition}(\textbf{Closure properties of fibrations})
		Let $P : \cat E \to \cat X$ be a fibration.
		\begin{enumerate}
			\item If $Q: \cat F \to \cat E$ is a fibration, then the composite functor $Q \comp P : \cat F \to \cat X$ is a fibration.
			\item If $Q: \cat F \to \cat X$ is any functor, then the pullback
			% https://q.uiver.app/#q=WzAsNCxbMCwwLCJRXipcXGNhdCBFIl0sWzAsMSwiXFxjYXQgRiJdLFsxLDEsIlxcY2F0IEIiXSxbMSwwLCJcXGNhdCBFIl0sWzAsMSwiUV4qUCIsMl0sWzEsMiwiUSJdLFszLDIsIlAiXSxbMCwzXSxbMCwyLCIiLDEseyJzdHlsZSI6eyJuYW1lIjoiY29ybmVyIn19XV0=
			\[\begin{tikzcd}[sep=4ex]
			{Q^*\cat E} & {\cat E} \\
			{\cat F} & {\cat X}
			\arrow["{Q^*P}"', from=1-1, to=2-1]
			\arrow["Q", from=2-1, to=2-2]
			\arrow["P", from=1-2, to=2-2]
			\arrow[from=1-1, to=1-2]
			\arrow["\lrcorner"{anchor=center, pos=0.125}, draw=none, from=1-1, to=2-2]
			\end{tikzcd}\]
			is a fibration.
		\end{enumerate}
	\end{proposition}

	Dualizing the notion of cartesian arrow gives rise to \emph{opcartesian arrows}, and dualizing the notion of fibration yields \emph{Grothendieck opfibrations}. In an opfibration $P: \cat E \to \cat X$, we can lift arrows $f : X \to Y$ given elements $U  :  \de P(X)$ to opcartesian arrows $U \to f_*U$ over $f$. While cartesian arrows satisfy a terminality property, opcartesian arrows satisfy an \emph{initiality} property.

	Like fibrations, opfibrations are closed under composition and pullback.

	\begin{example}(\textbf{Examples of fibrations})\label{ex:fibs}
		For a category $\cat C$, we denote by $\cat C^\downarrow$ the category of morphisms in $\cat C$ with arrows given by the commutative squares.
		\begin{enumerate}
			\item \textbf{Domain fibration}. For any category $\cat C$, the domain projection $\dom := \dom_{\cat C} : {\cat C}^\downarrow \to \cat C$ is a fibration. For $u:I \to J$ in $\cat X$, a cartesian lift of $u$ with respect to an arrow $g : J\to Y$ is given by the square:
			% https://q.uiver.app/#q=WzAsNCxbMCwwLCJJIl0sWzAsMSwiWCJdLFsxLDEsIlkiXSxbMSwwLCJKIl0sWzAsMSwidWciLDJdLFsxLDIsIiIsMCx7ImxldmVsIjoyLCJzdHlsZSI6eyJoZWFkIjp7Im5hbWUiOiJub25lIn19fV0sWzAsMywidSJdLFszLDIsImciXV0=
			\[\begin{tikzcd}[sep=4ex]
			I & J \\
			X & Y
			\arrow["u\comp g"', from=1-1, to=2-1]
			\arrow[Rightarrow, no head, from=2-1, to=2-2]
			\arrow["u", from=1-1, to=1-2]
			\arrow["g", from=1-2, to=2-2]
			\end{tikzcd}\]
			\item \textbf{Codomain fibration}. If $\cat C$ has all pullbacks, then the codomain projection $\cod := \cod_{\cat C} : \cat C^\downarrow \to \cat C$ is a fibration, also called the \emph{fundamental fibration} of $\cat C$. If $\cat C$ is locally cartesian closed, the codomain fibration gives rise to a model of \emph{dependent type theory} with $\Sigma$- and $\Pi$-types~\cite{seely1984locally}. The idea is to interpret a type $X$ in context (\emph{i.e.}, a list of free variables) $I$ as a morphism $X \to I$, hence an object of the slice $\cat C/I$, which is the fiber of $\cod$ at $I :  \cat C$.

			Given an arrow $u : J \to I$ in $\cat C$ together with a map $p : X \to I$, a $\cod$-cartesian lift of $u$ with respect to $p$ is given by pullback:
			\[\begin{tikzcd}[sep=4ex]
			{u^*X} & X \\
			J & I
			\arrow["u", from=2-1, to=2-2]
			\arrow["p", from=1-2, to=2-2]
			\arrow[from=1-1, to=2-1]
			\arrow[from=1-1, to=1-2]
			\arrow["\lrcorner"{anchor=center, pos=0.125}, draw=none, from=1-1, to=2-2]
			\end{tikzcd}\]
			This corresponds to substituting the variables in context $I$ along $u$.
			\item \textbf{Simple fibration}. Whereas the codomain fibration models arbitrary dependent types in context $I$, we might want to restrict to the case of \emph{simple types}. A \emph{constant type} is a dependent type over the terminal context $1$, hence a map $X \to 1$, or equivalently, just an object $X$. We can then view $X$ as a type over any context, by considering the pullback:
			\[\begin{tikzcd}[sep=4ex]
			{I \times X} & X \\
			I & 1
			\arrow[from=2-1, to=2-2]
			\arrow[from=1-2, to=2-2]
			\arrow[from=1-1, to=2-1]
			\arrow[from=1-1, to=1-2]
			\arrow["\lrcorner"{anchor=center, pos=0.125}, draw=none, from=1-1, to=2-2]
			\end{tikzcd}\]
			Therefore, as soon as $\cat C$ admits finite products, the \emph{simple fibration} is given by the subfunctor of $\cod_{\cat C}$ spanned by product projections $I \times X \to I$.
			It can also be described, up to fibered equivalence, as follows. Define the category $\Simple{\cat C}$ by taking as objects pairs $(I,X)$ of $\cat C$ and as morphisms $(I,X) \to (J,Y)$ pairs of morphisms $(u : I \to J,f : I \times X \to Y)$ in $\cat C$.
			The simple fibration $\simple : \Simple{\cat C} \to \cat C$ is the functor defined by $\simple(I,X) := I$ and $\simple(u,f) := u$. This is indeed a fibration, where a cartesian lift of $u : I \to J$ with respect to $(J,Y)$ is given by the map $(u,\pi_Y): (I,Y) \to (J,Y)$ with $\pi_Y : I \times Y \to Y$.
			\item \textbf{Subobject fibration}. Let $\cat C$ have all pullbacks. We can consider the full subfibration
			\begin{equation}
				\begin{tikzcd}[sep=4ex]
					{{\cat C}^\downarrow\vert_\mono} && {\cat C^\downarrow} \\
					& {\cat C}
					\arrow[hook, from=1-1, to=1-3]
					\arrow[from=1-1, to=2-2]
					\arrow[from=1-3, to=2-2]
				\end{tikzcd}
			\end{equation}
			of monomorphisms in $\cat C$. Quotienting out by the relation that two monomorphisms shall be equal if and only if they factor through each other (necessarily through a monomorphism), we obtain the \emph{subobject fibration} $\sub:\SubFib{\cat C} \to \cat C$. The subobject fibration is typically used as an elementary example of fibration of predicates on the objects of $\cat C$.
			\item \textbf{Family fibration}. Consider a category $\cat X$ with pullbacks and $P : \cat E \to \cat X$ a fibration.

			Let $\Fam(P) : P \downarrow \cat X \to \cat X$ be the functor constructed as follows:
			\begin{equation}
			\label{eq:fam}
				\begin{tikzcd}[sep=4ex]
					{P\downarrow \cat X} && {\cat E} \\
					{\cat X^\downarrow} && {\cat X} \\
					\cat X
					\arrow[from=1-1, to=2-1]
					\arrow["\dom", from=2-1, to=2-3]
					\arrow[from=1-1, to=1-3]
					\arrow["P", from=1-3, to=2-3]
					\arrow["\cod", from=2-1, to=3-1]
					\arrow["{\Fam(P)}"', curve={height=24pt}, from=1-1, to=3-1]
					\arrow["\lrcorner"{anchor=center, pos=0.125}, draw=none, from=1-1, to=2-3]
				\end{tikzcd}
			\end{equation}
			We call $\Fam(P) : P \downarrow \cat X \to \cat X$ the \emph{family fibration} of $P$. Indeed this is a fibration: any pullback of $P$ is a fibration, so is $\dom^*P$. Since $\cod$ is a fibration, as $\cat X$ has all pullbacks, then so is their composite $\cod \comp \dom^*P= \Fam(P)$. Objects in the fiber $\de\Fam(P)(I)$ are given by pairs $(f:X \to I, \alpha  :  \de P(X))$. In fact, if $P$ is an arbitrary functor, $\Fam(P)$ is always an \emph{opfibration}, the \emph{free opfibration} generated by $P$. Thus, $\Fam(P)$ is both a fibration and an opfibration. One can show that opcartesian transport serves as taking sums (or existential quantification), and that they commute with cartesian transport, serving as reindexing, in a canonical way (via the \emph{Beck--Chevalley isomorphisms}). Hence, the family construction on a fibration $P$ can be understood as the result of freely adding sums to $P$, which are \emph{internally indexed}, \emph{i.e.}, indexed by the base category $\cat X$. Accordingly, we also write $\Sum(P)$ for $\Fam(P)$.

			One could replace $\dom$ in the above definition by a full subfibration of itself, corresponding to a version of the family construction that is the sum completion only with respect to a chosen class of morphisms of $\cat X$.
		\end{enumerate}
	\end{example}

	\begin{definition}(\textbf{The vertical--cartesian factorization system, {\cite{myers_cartesian_2021}}})
		Let $P: \cat E \to \cat X$ be a fibration.
		The \textbf{vertical--cartesian factorization system} $(\ver, \car)$ induced by $P$ on $\cat E$ is an orthogonal factorization system where
		\begin{enumerate}
			\item the left class (\textbf{vertical morphisms}) is given by the maps projecting to isomorphisms and
			\item the right class (\textbf{cartesian morphisms}) is given by cartesian maps.
		\end{enumerate}
	\end{definition}

	In practice, this means every arrow of $\cat E$ can be factored in a unique way as a vertical map followed by a cartesian one:
	\begin{equation}
		% https://q.uiver.app?q=WzAsMyxbMCwwLCJcXGNkb3QiXSxbNCwwLCJcXGNkb3QiXSxbMiwwLCJcXGNkb3QiXSxbMCwxLCJcXHZhcnBoaSIsMCx7ImN1cnZlIjotNX1dLFswLDIsIlxcdmFycGhpXlxcbmF0IiwyXSxbMiwxLCJcXHZhcnBoaV5cXGZsYXQiLDJdXQ==
		\begin{tikzcd}[ampersand replacement=\&,sep=scriptsize]
			\cdot \&\& \cdot \&\& \cdot
			\arrow["\varphi", curve={height=-25pt}, dashed, from=1-1, to=1-5]
			\arrow["{\varphi^\ver}"', from=1-1, to=1-3]
			\arrow["{\varphi^\car}"', from=1-3, to=1-5]
		\end{tikzcd}
	\end{equation}

	This factorization is manifest when we present $P$ as the Grothendieck construction of its fibers: in that case, morphisms in the total category are literally pairs $\lens{f^\ver}{f^\car}$ where $f^\car$ is cartesian and $f^\ver$ is vertical.

	The most important property of this factorization system is the following:

	\begin{proposition}
		Let $P:\cat E \to \cat X$ be a fibration. In $\cat E$, every cospan:
		\begin{equation}
			% https://q.uiver.app/#q=WzAsMyxbMSwwLCJcXGNkb3QiXSxbMCwxLCJcXGNkb3QiXSxbMSwxLCJcXGNkb3QiXSxbMCwyLCJcXHZlciJdLFsxLDIsIlxcY2FyIiwyXV0=
			\begin{tikzcd}[ampersand replacement=\&]
				\& \cdot \\
				\cdot \& \cdot
				\arrow["\ver", from=1-2, to=2-2]
				\arrow["\car"', from=2-1, to=2-2]
			\end{tikzcd}
		\end{equation}
		admits a pullback:
		\begin{equation}
		\label{eq:carversq}
			% https://q.uiver.app/#q=WzAsNCxbMSwwLCJcXGNkb3QiXSxbMCwxLCJcXGNkb3QiXSxbMSwxLCJcXGNkb3QiXSxbMCwwLCJcXGNkb3QiXSxbMCwyLCJcXHZlciJdLFsxLDIsIlxcY2FyIiwyXSxbMywxLCJcXHZlciIsMl0sWzMsMCwiXFxjYXIiXSxbMywyLCIiLDEseyJzdHlsZSI6eyJuYW1lIjoiY29ybmVyIn19XV0=
			\begin{tikzcd}[ampersand replacement=\&]
				\cdot \& \cdot \\
				\cdot \& \cdot
				\arrow["\ver", from=1-2, to=2-2]
				\arrow["\car"', from=2-1, to=2-2]
				\arrow["\ver"', from=1-1, to=2-1]
				\arrow["\car", from=1-1, to=1-2]
				\arrow["\lrcorner"{anchor=center, pos=0.125}, draw=none, from=1-1, to=2-2]
			\end{tikzcd}
		\end{equation}
		Conversely, every commutative square like~\eqref{eq:carversq} is a pullback square.
	\end{proposition}

	\begin{construction}
		Given a fibration $P : \cat E \to \cat X$, its \textbf{dual} (or \textbf{fiberwise opposite}) is another fibration $\dual P : \dual{\cat E} \to \cat X$ constructed as follows.

		The category $\dual{\cat E}$ is the category that has the same objects as $\cat E$ and as morphisms equivalence classes of vertical-cartesian spans in $\cat E$, i.e.~spans in $\cat E$ of the form
		\begin{equation}
			% https://q.uiver.app/?q=WzAsMyxbMCwxLCJYIl0sWzEsMCwiXFxjZG90ICJdLFsyLDEsIlkiXSxbMSwyLCJcXGNhciJdLFsxLDAsIlxcdmVyIiwyXV0=
			\begin{tikzcd}[ampersand replacement=\&, sep=small]
				\& {\cdot } \\
				\cdot \&\& \cdot
				\arrow["\car", from=1-2, to=2-3]
				\arrow["\ver"', from=1-2, to=2-1]
			\end{tikzcd}
		\end{equation}
		They compose because vertical arrows can always be pulled back against cartesian arrows:
		\begin{equation}
			% https://q.uiver.app/#q=WzAsNixbMCwyLCJYIl0sWzEsMSwiXFxjZG90ICJdLFsyLDIsIlkiXSxbMywxLCJcXGNkb3QiXSxbNCwyLCJaIl0sWzIsMCwiXFxjZG90Il0sWzEsMiwiXFxjYXIiLDAseyJsYWJlbF9wb3NpdGlvbiI6ODB9XSxbMSwwLCJcXHZlciIsMl0sWzMsMiwiXFx2ZXIiLDIseyJsYWJlbF9wb3NpdGlvbiI6MzB9XSxbMyw0LCJcXGNhciJdLFs1LDEsIlxcdmVyIiwyXSxbNSwzLCJcXGNhciJdLFs1LDIsIiIsMix7InN0eWxlIjp7Im5hbWUiOiJjb3JuZXIifX1dXQ==
		\begin{tikzcd}[ampersand replacement=\&, sep=small]
			\&\& \cdot \\
			\& {\cdot } \&\&\cdot \\
			X \&\& Y \&\& Z
			\arrow["\car"{pos=0.8}, from=2-2, to=3-3]
			\arrow["\ver"', from=2-2, to=3-1]
			\arrow["\ver"'{pos=0.3}, from=2-4, to=3-3]
			\arrow["\car", from=2-4, to=3-5]
			\arrow["\ver"', from=1-3, to=2-2]
			\arrow["\car", from=1-3, to=2-4]
			\arrow["\lrcorner"{anchor=center, pos=0.125, rotate=-45}, draw=none, from=1-3, to=3-3]
		\end{tikzcd}
		\end{equation}
		Hence this defines a category, and $P$ still works as a projection onto $\cat X$.

		The fibration structure comes from noticing one can still lift arrows to cartesian arrows. More easily, since our fibrations are cloven, we can define
		\begin{equation}
			\dual P := \int (\de P \comp (-)^\op)
		\end{equation}
		which is a fibration by construction.
		This point of view will be also useful later on.
	\end{construction}

	Alternatively, one can define the arrows $E \to D$ in the dual fibration as pairs
	\[(f : P(D) \to P(E), \psi : X \to f^*E),\]
	where $\psi$ is vertical. From this construction one can readily see that the ensuing total category $\cat E ^\vee$ is locally small if $\cat X$ and $\cat E$ are.

	The dual fibration construction extends to a functor $\dual{(-)} : \Fib(\cat X) \to \Fib(\cat X)$.

	Recall that a map of fibrations $F : Q \to P$ is a functor between the total spaces of $P$ and $Q$ that makes the obvious triangle commute and that respects the cartesian-vertical factorization system:
	\begin{equation}
		% https://q.uiver.app/?q=WzAsMyxbMiwwLCJcXEUiXSxbMSwxLCJcXEIiXSxbMCwwLCJcXEQiXSxbMiwwLCJmIl0sWzIsMSwicSIsMl0sWzAsMSwicCJdXQ==
		\begin{tikzcd}[ampersand replacement=\&, sep=small]
			\cat D \&\& \cat E \\
			\& \cat X
			\arrow["F", from=1-1, to=1-3]
			\arrow["Q"', from=1-1, to=2-2]
			\arrow["P", from=1-3, to=2-2]
		\end{tikzcd}
	\end{equation}
	Explicitly, it means $F$ sends cartesian arrows to cartesian arrows.
	This is equivalent to say $F$ defines functors between the fibers of $P$ and $Q$: for a given $X  :  \cat X$, $\de F(X) : \de P(X) \to \de Q(X)$ is a well-defined functor mapping $P$-vertical arrows over $X$ to $Q$-vertical arrows over $X$, and such that for every arrow $f : X' \to X$ in $\cat X$, the following square commutes:
	\begin{equation}
		% https://q.uiver.app?q=WzAsNCxbMCwwLCJcXGRlIFAoQikiXSxbMSwwLCJcXGRlIFEoQikiXSxbMCwxLCJcXGRlIFAoQicpIl0sWzEsMSwiXFxkZSBRKEInKSJdLFswLDEsIlxcZGUgRl9CIl0sWzIsMywiXFxkZSBGX3tCJ30iXSxbMCwyLCJcXGRlIFAoZikiLDJdLFsxLDMsIlxcZGUgUShmKSJdXQ==
		\begin{tikzcd}[ampersand replacement=\&]
			{\de P(X)} \& {\de Q(X)} \\
			{\de P(X')} \& {\de Q(X')}
			\arrow["{\de F(X)}", from=1-1, to=1-2]
			\arrow["{\de F(X')}", from=2-1, to=2-2]
			\arrow["{\de P(f)}"', from=1-1, to=2-1]
			\arrow["{\de Q(f)}", from=1-2, to=2-2]
		\end{tikzcd}
	\end{equation}

	Now if we take the fiberwise opposites of $P$ and $Q$, we can still define $\dual F$ by applying $F$ to both legs of the spans. Since $F$ respects the cartesian factorization systems of the fibrations, the image span is still of the required form.

	\section{Lens, optics and Dialectica categories}
	\label{sec:examples}
	\subsection{Lenses}
	In their most essential form, lenses are just morphisms in a dualized fibration.
	Originally, the name lenses referred to morphisms in the dual of the simple fibration, so that $\Lens(\cat C) := \dual{\Simple{\cat C}}$, where $\cat C$ is a category with finite products.
	Thus a \emph{(simple) lens} looks like this:
	\begin{equation}\label{eq:simple.lens}
		% https://q.uiver.app?q=WzAsNixbMSwxLCJVIl0sWzIsMSwiViJdLFsyLDAsIlYgXFx0aW1lcyBZIl0sWzAsMSwiVSJdLFswLDAsIlUgXFx0aW1lcyBYIl0sWzEsMCwiVSBcXHRpbWVzIFkiXSxbMCwxLCJmIiwwLHsiY29sb3VyIjpbMCw2MCw2MF19LFswLDYwLDYwLDFdXSxbMiwxLCJcXHBpX1YiLDJdLFszLDAsIiIsMCx7ImxldmVsIjoyLCJzdHlsZSI6eyJoZWFkIjp7Im5hbWUiOiJub25lIn19fV0sWzUsMiwiZiBcXHRpbWVzIFkiXSxbNSwxLCIiLDAseyJzdHlsZSI6eyJuYW1lIjoiY29ybmVyIn19XSxbNSwwLCJcXHBpX1UiLDJdLFs0LDMsIlxccGlfVSIsMl0sWzUsNCwiZl5cXHNoYXJwIiwyLHsiY29sb3VyIjpbMCw2MCw2MF19LFswLDYwLDYwLDFdXV0=
		\begin{tikzcd}[ampersand replacement=\&]
			{U \times X} \& {U \times Y} \& {V \times Y} \\
			U \& U \& V
			\arrow["f", dashed, from=2-2, to=2-3]
			\arrow["{\pi_V}"', from=1-3, to=2-3]
			\arrow[Rightarrow, no head, from=2-1, to=2-2]
			\arrow["{f \times Y}", from=1-2, to=1-3]
			\arrow["\lrcorner"{anchor=center, pos=0.125}, draw=none, from=1-2, to=2-3]
			\arrow["{\pi_U}"', from=1-2, to=2-2]
			\arrow["{\pi_U}"', from=1-1, to=2-1]
			\arrow["{f^\sharp}"', dashed, from=1-2, to=1-1]
		\end{tikzcd}
	\end{equation}
	In~\cite{spivak_generalized_2019}, Spivak shows the name lenses deserves to be applied to much more general morphisms
	Chiefly, \emph{dependent lenses} are morphisms in the dual of the codomain fibration of a category with pullbacks:
	\begin{equation}
		% https://q.uiver.app?q=WzAsNixbMSwxLCJVIl0sWzIsMSwiViJdLFsyLDAsIlkiXSxbMCwxLCJVIl0sWzAsMCwiWCJdLFsxLDAsIlUgXFx0aW1lc19WIFkiXSxbMCwxLCJmIiwwLHsiY29sb3VyIjpbMCw2MCw2MF19LFswLDYwLDYwLDFdXSxbMiwxXSxbMywwLCIiLDAseyJsZXZlbCI6Miwic3R5bGUiOnsiaGVhZCI6eyJuYW1lIjoibm9uZSJ9fX1dLFs1LDIsImZfWSJdLFs1LDEsIiIsMCx7InN0eWxlIjp7Im5hbWUiOiJjb3JuZXIifX1dLFs1LDBdLFs0LDNdLFs1LDQsImZeXFxzaGFycCIsMix7ImNvbG91ciI6WzAsNjAsNjBdfSxbMCw2MCw2MCwxXV1d
		\begin{tikzcd}[ampersand replacement=\&]
			X \& {U \times_V Y} \& Y \\
			U \& U \& V
			\arrow["f", dashed, from=2-2, to=2-3]
			\arrow[from=1-3, to=2-3]
			\arrow[Rightarrow, no head, from=2-1, to=2-2]
			\arrow["{f_Y}", from=1-2, to=1-3]
			\arrow["\lrcorner"{anchor=center, pos=0.125}, draw=none, from=1-2, to=2-3]
			\arrow[from=1-2, to=2-2]
			\arrow[from=1-1, to=2-1]
			\arrow["{f^\sharp}"', dashed, from=1-2, to=1-1]
		\end{tikzcd}
	\end{equation}

	\begin{definition}(\textbf{$P$-lenses})
	\label{def:lens}
		Let $P: \cat E \to \cat X$ be a fibration.
		Its dual fibration $\dual P: \dual{\cat E} \to \cat X$ is called the \textbf{fibration of $P$-lenses}, and morphisms in $\dual{\cat E}$ are \textbf{$P$-lenses}.
	\end{definition}

	\subsection{Dialectica categories}
	Dialectica categories~\cite{de_paiva_dialectica_1989} are made of morphisms very similar to lenses, and also have a formal essence which is very simple to capture with very little structure. For the categories we are presenting in due course, we are employing type-theoretic notation.
	\begin{definition}
	\label{def:dc}
		The \textbf{(simple) Dialectica category} $\sDial(\cat C)$ has
		\begin{enumerate}
			\item as objects, triples
			\begin{equation}
				(U : \cat C, X : \cat C, \alpha : \Sub(U \times X))
			\end{equation}
			meaning $\alpha$ is a subobject of $U \times X$,
			\item as morphisms $(U,X,\alpha) \to (V,Y,\beta)$, triples
			\begin{eqalign}
				&f: U \to V,\\
				&f^\sharp : U \times Y \to X,\\
				\forall u:U, y:Y, \quad &\alpha(u, f^\sharp(u, y)) \subseteq \beta(f(u), y)
			\end{eqalign}
		\end{enumerate}
	\end{definition}

	Hence $\sDial(\cat C)$ is a category of lenses `augmented with predicates'.

	\subsection{Optics}
	Optics are a generalization of lenses-as-data accessors introduced by the functional programming community to treat data structures other than record types (see~\cite{pickering_profunctor_2017}).

	Recall from~\eqref{eq:simple.lens} that a simple lens $\lens{f^\sharp}{f} : \lens{X}{U} \opticto \lens{Y}{V}$ is comprised of a \emph{forward map} $f:U \to V$ and a \emph{backward map} $f^\sharp : U \times Y \to X$.
	The extra argument $U$ in this latter map is called \emph{residual}, and in optics its role becomes central: turned into a first-class object, it is explicitly `written' in the forward part and `read' in the backward part.

	For the sake of brevity, we will only sketch the relevant definitions, pointing out the relevant precise references.

	Optics are defined via the concept of \emph{action of a monoidal category} (also known as \emph{actegories},~\cite{capucci_actegories_2023}):

	\begin{definition}(\textbf{Action of a monoidal category})
	\label{def:actegories}
		Let $(\cat M, I, \otimes)$ be a monoidal category.
		A \emph{(left) $\cat M$-action on the category $\cat C$} is a functor $\bullet : \cat M \times \cat C \to \cat C$ together with coherent isomorphisms witnessing $I \bullet A \cong A$ and $M \bullet (N \bullet A) \cong (M \otimes N) \bullet A$.
		The action is \emph{strict} when such isomorphisms are in fact identities.
	\end{definition}

	Like for monoidal categories, we refer to actegories by synecdoche using just the name of the carrier category.
	A full definition and a comprehensive account of actegories can be found in~\cite{capucci_actegories_2023}.
	Note that all actegories can be stricified.
	We work with those for simplicity.

	\newcommand{\PreOptic}{\cat{PreOptic}}
	\begin{definition}(\textbf{Preoptic})
	\label{def:preoptics}
		Fix a strict monoidal category $\cat M$ and two strict left $\cat M$-actions $(\cat C,\bullet)$ and $(\cat D, \circ)$.
		Let $A,S : \cat C$ and $B,T  :  \cat D$ be objects.
		A \emph{preoptic} $\lens{B}{A} \opticto \lens{T}{S}$ is a triple $(M,f,f^\sharp)$ where $M  : \cat M$ is called the \emph{residual}, $f:A \to M \bullet S$ is a morphism in $\cat C$ called the \emph{view}, and $f^\sharp : M \circ T \to B$ is a morphism in $\cat D$ called the \emph{update}.
		The category of \emph{preoptics} $\PreOptic_{\bullet,\circ}$ has objects pairs of an object of $\cat C$ and an object of $\cat D$, and morphisms preoptics.
	\end{definition}

	Optics are the category obtained by quotienting the hom-categories of $\PreOptic$ by an equivalence relation that allows morphisms to `slide' across the residual, roughly enforcing the fact that computing something from some inputs and then saving it in memory is equivalent to saving in memory the inputs and then later computing from those.

	By varying the different pieces of data, optics can express a plethora of data accessing patterns. For a survey and more detailed definitions, we direct the reader to~\cite{clarke_profunctor_2024}.

	\begin{remark}
		Using preoptics instead of optics is a convenient device to avoid working with \emph{2-optics} instead, in which the quotient hinted at above is replaced by explicit witnesses.
		Avoiding the quotient, however, seems necessary for the fibrational manipulation of optics, see e.g.~\cite{braithwaite_fibre_2021}.
	\end{remark}

	\section{Dialenses as a common framework}
	\label{sec:dialenses}

	\subsection{Towers of fibrations and their duals}
	The data needed to construct a Dialectica category is easily seen to amount to a tower of fibrations:
	\begin{equation}
		% https://q.uiver.app?q=WzAsNixbMCwwLCJcXGNhdCBQIl0sWzAsMSwiXFxjYXQgRSJdLFswLDIsIlxcY2F0IEMiXSxbMSwyLCJmOlUgXFx0byBWIl0sWzEsMSwiZl5cXHNoYXJwOiBmXipZIFxcdG8gWCJdLFsxLDAsImZeXFxkb3VibGVzaGFycCA6IChmXlxcc2hhcnApXipcXGFscGhhIFxcdG8gZl4qXFxiZXRhIl0sWzEsMiwicCIsMl0sWzAsMSwicSIsMl1d
		\begin{tikzcd}[ampersand replacement=\&]
			{\cat P} \& {`\alpha \subseteq \beta\text{'}} \\
			{\cat E} \& {f^\sharp: f^*Y \to X} \\
			{\cat C} \& {f:U \to V}
			\arrow["P"', from=2-1, to=3-1]
			\arrow["Q"', from=1-1, to=2-1]
		\end{tikzcd}
	\end{equation}
	Thus we can wonder if somehow the simple construction producing $P$-lenses from $P$ extends to the extra data of $Q$.

	First, notice such a tower of fibrations can be seen as exhibiting $Q$ as a morphism of fibrations over $\cat C$:
	\begin{equation}
		% https://q.uiver.app?q=WzAsMyxbMCwwLCJcXGNhdCBQIl0sWzIsMCwiXFxjYXQgRSJdLFsxLDEsIlxcY2F0IEMiXSxbMSwyLCJwIl0sWzAsMSwicSJdLFswLDIsInEgXFxjb21wIHAiLDJdXQ==
		\begin{tikzcd}[ampersand replacement=\&]
			{\cat P} \&\& {\cat E} \\
			\& {\cat C}
			\arrow["P", from=1-3, to=2-2]
			\arrow["Q", from=1-1, to=1-3]
			\arrow["{Q \comp P}"', from=1-1, to=2-2]
		\end{tikzcd}
	\end{equation}
	In fact the triangle trivially commutes and $Q$ respects cartesian arrows since the cartesian arrows of $\cat P$ \emph{as a category fibered over $\cat C$} are exactly the cartesian arrows of $\cat P$ (as a category fibered over $\cat E$) which are over cartesian arrows of $\cat E$.

	\begin{remark}
		Something stronger is true: $Q$ is actually a \emph{fibration of fibrations}, i.e.~a `fiberwise fibration', since B\'enabou showed that fibrations in $\Fib(\cat C)$ are given by maps of fibrations \emph{which are themselves fibrations}.
		A published account of this fact can be found in~\cite{hermida_properties_1999}.
	\end{remark}

	Then, since $\dual{(-)}$ is a functorial construction $\Fib(\cat C) \to \Fib(\cat C)$, it can be applied to the whole triangle:
	\begin{equation}
		% https://q.uiver.app?q=WzAsMyxbMCwwLCJcXGNhdCBQIl0sWzIsMCwiXFxjYXQgRSJdLFsxLDEsIlxcY2F0IEMiXSxbMSwyLCJwIl0sWzAsMSwicSJdLFswLDIsInEgXFxjb21wIHAiLDJdXQ==
		\begin{tikzcd}[ampersand replacement=\&]
			{\dual{\cat P}} \&\& {\dual{\cat E}} \\
			\& {\cat C}
			\arrow["\dual P", from=1-3, to=2-2]
			\arrow["\dual Q", from=1-1, to=1-3]
			\arrow["\dual {(Q \comp P)}"', from=1-1, to=2-2]
		\end{tikzcd}
	\end{equation}

	To understand how $\dual{\cat P}$ looks like now, we need to stop for a moment and understand the factorization systems at play in this situation.

	Remember $P$ induces a factorization system on $\cat E$ which $\dual{\cat E}$ inherits as ($P$-vertical$^\op$, $P$-cartesian) (also denoted as $(\ver_P^\op, \car_P)$).
	On $\cat P$ we have a more refined factorization system since we have three kinds of arrows:
	\begin{enumerate}
		\item \textbf{$Q$-cartesian arrows}, which are cartesian lifts of $\cat E$-arrows, and come in two subcategories:
		\begin{enumerate}
			\item \textbf{$P$-cartesian arrows}, which are cartesian lifts of $\cat C$-arrows, hence cartesian $\cat E$-arrows,
			\item \textbf{$P$-vertical arrows}, which are cartesian lifts of vertical $\cat E$-arrows,
		\end{enumerate}
		\item \textbf{$Q$-vertical arrows}, which are in the fibers of $P$.
	\end{enumerate}

	This forms a ternary factorization system \textbf{($Q$-vertical, $P$-vertical, $P$-cartesian)}, or $(\ver_Q, \ver_P, \car_P)$, which means every arrow $\varphi$ in $\cat P$ factors in three parts, uniquely up to unique isomorphism:
	\begin{equation}
		% https://q.uiver.app?q=WzAsNCxbNCwwLCJcXGNkb3QiXSxbMiwwLCJcXGNkb3QiXSxbNiwwLCJcXGNkb3QiXSxbMCwwLCJcXGNkb3QiXSxbMywxLCJcXHZhcnBoaV57XFx2ZXJfUX0iLDJdLFsxLDAsIlxcdmFycGhpXntcXHZlcl9QfSIsMl0sWzAsMiwiXFx2YXJwaGlee1xcY2FyX1B9IiwyXSxbMywyLCJcXHZhcnBoaSIsMSx7ImN1cnZlIjotNSwic3R5bGUiOnsiYm9keSI6eyJuYW1lIjoiZGFzaGVkIn19fV0sWzEsMiwiXFx2YXJwaGlee1xcY2FyX1F9IiwxLHsiY3VydmUiOi0zLCJzdHlsZSI6eyJib2R5Ijp7Im5hbWUiOiJkYXNoZWQifX19XV0=
		\begin{tikzcd}[ampersand replacement=\&,sep=scriptsize]
			\cdot \&\& \cdot \&\& \cdot \&\& \cdot
			\arrow["{\varphi^{\ver_Q}}"', from=1-1, to=1-3]
			\arrow["{\varphi^{\ver_P}}"', from=1-3, to=1-5]
			\arrow["{\varphi^{\car_P}}"', from=1-5, to=1-7]
			\arrow["\varphi"{description}, curve={height=-50pt}, dashed, from=1-1, to=1-7]
			\arrow["{\varphi^{\car_Q}}"{description}, curve={height=-25pt}, dashed, from=1-3, to=1-7]
		\end{tikzcd}
	\end{equation}
	When we turn around the fibers of $\cat E$, we end up swapping $Q$-vertical and $P$-vertical arrows:
	\begin{equation}
		% https://q.uiver.app/#q=WzAsNCxbNCwwLCJcXGNkb3QiXSxbMiwwLCJcXGNkb3QiXSxbNiwwLCJcXGNkb3QiXSxbMCwwLCJcXGNkb3QiXSxbMSwzLCJcXHZhcnBoaV57XFx2ZXJfUX0iXSxbMCwxLCJcXHZhcnBoaV57XFx2ZXJfUH0iXSxbMCwyLCJcXHZhcnBoaV57XFxjYXJfUH0iLDJdLFszLDIsIlxcdmFycGhpIiwxLHsiY3VydmUiOi01fV1d
		\begin{tikzcd}[ampersand replacement=\&,sep=scriptsize]
			\cdot \&\& \cdot \&\& \cdot \&\& \cdot
			\arrow["{\varphi^{\ver_Q}}", from=1-3, to=1-1]
			\arrow["{\varphi^{\ver_P}}", from=1-5, to=1-3]
			\arrow["{\varphi^{\car_P}}"', from=1-5, to=1-7]
			\arrow["\varphi"{description}, curve={height=-30pt}, dashed, from=1-1, to=1-7]
		\end{tikzcd}
	\end{equation}
	Hence, on $\dual{\cat P}$ we end up with a ternary factorization system ($P$-vertical$^\op$, $Q$-vertical$^\op$, $P$-cartesian), i.e. where $P$- and $Q$-vertical arrows are swapped (both as arrows and as classes in the ternary factorization system).

	We can understand this factorization system as an \textbf{ambifibration}\footnote{The definition and name of ambifibration seems to be due to Kock and Joyal, but unpublished. See~\cite{kock_comment_2010}.} structure on $\dual Q$:

	\begin{definition}
		Let $(L, R)$ be a factorization system on $\cat E$. An \textbf{ambifibration} $A : \cat P \to \cat E$ is a functor such that
		\begin{enumerate}
			\item every arrow in $L$ has an opcartesian lift ($A$ restrics to an opfibration on $L$)
			\item every arrow in $R$ has a cartesian lift ($A$ restrics to a fibration on $R$)
		\end{enumerate}
	\end{definition}
	This induces the ternary factorization system (opcartesian, vertical, cartesian) on $\cat P$:
	\begin{equation}
		% file:///home/jsb20179/data/software/quiver/src/index.html?q=WzAsMTAsWzMsMSwiQiJdLFs0LDEsIkMiXSxbMSwxLCJBIl0sWzIsMSwiQiJdLFs0LDAsIkMnIl0sWzEsMCwiQSciXSxbMiwwLCJcXGVsbF8qIEEnIl0sWzMsMCwicl4qQyciXSxbMCwwLCJcXGNhdCBGIl0sWzAsMSwiXFxjYXQgRCJdLFswLDEsIlxcaW4gUiIsMl0sWzIsMywiXFxpbiBMIiwyXSxbMywwLCIiLDIseyJsZXZlbCI6Miwic3R5bGUiOnsiaGVhZCI6eyJuYW1lIjoibm9uZSJ9fX1dLFs1LDYsIlxcb3BjYXIiLDJdLFs3LDQsIlxcY2FyIiwyXSxbNiw3LCJcXGluIFYiLDJdLFs4LDksImEiLDJdXQ==
		\begin{tikzcd}[ampersand replacement=\&]
			{\cat P} \& {A'} \& {\ell_* A'} \& {r^*C'} \& {C'} \\
			{\cat E} \& A \& B \& B \& C
			\arrow["{r \in R}"', from=2-4, to=2-5]
			\arrow["{\ell \in L}"', from=2-2, to=2-3]
			\arrow[Rightarrow, no head, from=2-3, to=2-4]
			\arrow["\opcar"', from=1-2, to=1-3]
			\arrow["\car"', from=1-4, to=1-5]
			\arrow["\ver"', from=1-3, to=1-4]
			\arrow["A"', from=1-1, to=2-1]
		\end{tikzcd}
	\end{equation}

	Looking at the tower of fibrations $\dual{\cat P} \nto{\dual Q} \dual{\cat E} \nto{\dual P} \cat C$, we see:
	\begin{enumerate}
		\item $\dual{\cat E}$ has the factorization system ($P$-vertical$^\op$, $P$-cartesian),
		\item $\dual Q$ is a fibration on the cartesian class of $\cat E$, given by $Q$, and became an opfibration on the vertical class since it acts like $\dual{Q}$ there.
		Hence $\dual{Q}$ is an ambifibration over $\dual{\cat E}$ with respect to the factorization system $(\ver^\op, \car)$.
	\end{enumerate}

	This is very close to what we are looking for: the data of a morphism in Dialectica is indeed that of a triple of morphisms, except the last one (first in the factorization) needs to point in the same direction as the first one (last in the factorization).
	This suggests a last dualization is needed, involving only the fibration $Q$.

	At this point, it's a good idea to introduce some new notation and terminology.
	First of all, let us make the notation $\dual{P} : \dual{\cat E} \to \cat X$ more precise by specifying the base with respect to which we are taking the dual:
	\begin{equation}
		\dual[\cat X]{P} : \dual[\cat X]{\cat E} \to \cat X.
	\end{equation}
	Moreover, we introduce the following construction:

	\begin{definition}
		Let $\tower{P} = (\cat E_n \nlongto{P_n} \cdots \nlongto{P_3} \cat E_2 \nlongto{P_2} \cat E_1 \nlongto{P_1} \cat E_0)$ be a sequence of $n$ fibrations (a  \textbf{tower of fibrations of height} $n$).
		The \textbf{iterated dual} of $\tower{P}$ is defined inductively on its length $1 \leq k < n$ as follows:
		\begin{itemize}
			\item $k=1:$
			\begin{equation}
				\dual[1]{(\cat E_1 \nlongto{P_1} \cat E_0)} = \dual[\cat E_0]{\cat E_1} \nlongto{\dual[\cat E_0]{P_1}} \cat E_0
			\end{equation}
			\item $1 < k < n:$
			\begin{equation}
				\dual[k+1]{(\cat E_{k+1} \nlongto{P_{k+1}} \cdots \nlongto{P_2} \cat E_1 \nlongto{P_1} \cat E_0)}
				=
				\dual[\cat E_0]{\left(\dual[k]{(\cat E_{k+1} \nlongto{P_{k+1}} \cdots \nlongto{P_{2}} \cat E_1)}\right)}\nlongto{\dual[\cat E_0]{P_1}} \cat E_0
			\end{equation}
		\end{itemize}
	\end{definition}

	There are many natural examples of towers of fibrations:

	\begin{example}
		If $\cat C$ is a cartesian monoidal category, then the total category of its associated simple fibration $\simple:\Simple{\cat C} \to \cat C$ is again cartesian, meaning we can form towers of arbitrary height by iterating the construction of simple fibrations.
	\end{example}

	\begin{example}
		Likewise, when $\cat C$ is finitely complete we can form towers $\cat C^{\downarrow^n} \to \cdots \to \cat C^\downarrow \to \cat C$ of arbitrary height where the $i$-th level is the category of \emph{$i$-dimensional commutative cubes} in $\cat C$.
	\end{example}

	\begin{example}
		More generally, suppose $P:\cat E \to \cat C$ is a \emph{full comprehension category}~\cite{jacobs_comprehension_1993}, meaning it admits a right adjoint $1:\cat C \to \cat E$ (a \emph{fibered terminal object}) as well as a further right adjoint (\emph{comprehension}) $\{-\}:\cat E \to \cat C$.
		By repeatedly pulling back $P$ along $\{-\}$ we can form towers of fibrations of arbitrary height.
	\end{example}

	\subsection{Dialenses}
	Let's see the result of an iterated dual of length $2$, i.e. $\dual[2]{(\cat E_2 \nto{P_2} \cat E_1 \nto{P_1} \cat E_0)}$.
	Below are unpacked the various bits the morphisms in $\dual[\cat E_1]{\dual[\cat E_0]{\cat E_2}}$ are made of: as you can see we end up with a ternary factorization corresponding to three changes in directions as we go up the sequence of fibrations:
	\begin{equation}
		% file:///home/jsb20179/data/software/quiver/src/index.html?q=WzAsMTUsWzMsNCwiVSJdLFs0LDQsIlYiXSxbMCw0LCJcXGNhdCBDIl0sWzAsMiwiXFxkdWFsW1xcY2F0IENde1xcY2F0IEVfMX0iXSxbNCwyLCJZIl0sWzMsMiwiZl4qWSJdLFsxLDQsIlUiXSxbMiw0LCJVIl0sWzIsMiwiZl4qWSJdLFsxLDIsIlgiXSxbMCwwLCJcXGR1YWxbXFxjYXQgRV8xXXtcXGR1YWxbXFxjYXQgQ117XFxjYXQgRV8yfX0iXSxbMSwwLCJcXGFscGhhIl0sWzQsMCwiXFxiZXRhIl0sWzIsMCwiKGZeXFxzaGFycCleKlxcYWxwaGEiXSxbMywwLCJmXipcXGJldGEiXSxbMCwxLCJmIiwwLHsiY29sb3VyIjpbMCw2MCw2MF19LFswLDYwLDYwLDFdXSxbNSw0LCJmX1kiXSxbNiw3LCIiLDAseyJsZXZlbCI6Miwic3R5bGUiOnsiaGVhZCI6eyJuYW1lIjoibm9uZSJ9fX1dLFs3LDAsIiIsMCx7ImxldmVsIjoyLCJzdHlsZSI6eyJoZWFkIjp7Im5hbWUiOiJub25lIn19fV0sWzgsOSwiZl5cXHNoYXJwIiwyLHsiY29sb3VyIjpbMCw2MCw2MF19LFswLDYwLDYwLDFdXSxbOCw1LCIiLDAseyJsZXZlbCI6Miwic3R5bGUiOnsiaGVhZCI6eyJuYW1lIjoibm9uZSJ9fX1dLFszLDIsIlxcZHVhbFtcXGNhdCBDXXtwXzF9IiwyXSxbMTAsMywiXFxkdWFsW1xcY2F0IEVfMV17XFxkdWFsW1xcY2F0IENde3BfMn19IiwyXSxbMTQsMTIsImZfXFxiZXRhIl0sWzEzLDExLCJmXlxcc2hhcnBfXFxhbHBoYSIsMl0sWzEzLDE0LCJmXlxcZG91Ymxlc2hhcnAiLDAseyJjb2xvdXIiOlswLDYwLDYwXX0sWzAsNjAsNjAsMV1dXQ==
		\begin{tikzcd}[ampersand replacement=\&, sep=scriptsize]
			{\dual[\cat E_1]{\dual[\cat E_0]{\cat E_2}}} \& \alpha \& {(f^\sharp)^*\alpha} \& {f^*\beta} \& \beta \\
			\\
			{\dual[\cat E_0]{\cat E_1}} \& X \& {f^*Y} \& {f^*Y} \& Y \\
			\\
			{\cat E_0} \& U \& U \& U \& V
			\arrow["f", color={red}, from=5-4, to=5-5]
			\arrow["{f_Y}", from=3-4, to=3-5]
			\arrow[Rightarrow, no head, from=5-2, to=5-3]
			\arrow[Rightarrow, no head, from=5-3, to=5-4]
			\arrow["{f^\sharp}"', color={red}, from=3-3, to=3-2]
			\arrow[Rightarrow, no head, from=3-3, to=3-4]
			\arrow["{\dual[\cat E_0]{P_1}}"', from=3-1, to=5-1]
			\arrow["{\dual[\cat E_1]{\dual[\cat E_0]{P_2}}}"', from=1-1, to=3-1]
			\arrow["{f_\beta}", from=1-4, to=1-5]
			\arrow["{f^\sharp_\alpha}"', from=1-3, to=1-2]
			\arrow["{f^\doublesharp}", color={red}, from=1-3, to=1-4]
		\end{tikzcd}
	\end{equation}

	To see this reproduces Definition~\ref{def:dc}, let's unpack this situation when $P_2 = \simple:\Simple{\cat C} \to \cat C$ and $P_1$ is the fibration of monos of the domain over $\Simple{\cat C}$, that is, the fiber over $\pi_U : U \times X \to U$ is given by the subobjects of $U \times X$:
	\begin{equation}
		% file:///home/jsb20179/data/software/quiver/src/index.html?q=WzAsMTUsWzIsMiwiWCJdLFszLDIsIlgiXSxbNCwyLCJYIixbMCw2MCw2MCwxXV0sWzUsMiwiWSIsWzAsNjAsNjAsMV1dLFs1LDEsIlYgXFx0aW1lcyBZIl0sWzQsMSwiViBcXHRpbWVzIFgiXSxbNSwwLCJcXGJldGEiXSxbMiwwLCJcXGFscGhhIl0sWzIsMSwiVSBcXHRpbWVzIFgiLFswLDYwLDYwLDFdXSxbMywxLCJWIFxcdGltZXMgWCIsWzAsNjAsNjAsMV1dLFszLDAsIihmXlxcc2hhcnApXipcXGFscGhhIixbMCw2MCw2MCwxXV0sWzQsMCwiZl4qXFxiZXRhIixbMCw2MCw2MCwxXV0sWzAsMiwiXFxjYXQgQyJdLFswLDEsIlxcY2F0IENeXFxkb3duYXJyb3dfXFxwcm9qIl0sWzAsMCwiXFxjYXQgUCJdLFswLDEsIiIsMix7ImxldmVsIjoyLCJzdHlsZSI6eyJoZWFkIjp7Im5hbWUiOiJub25lIn19fV0sWzEsMiwiIiwyLHsibGV2ZWwiOjIsInN0eWxlIjp7ImhlYWQiOnsibmFtZSI6Im5vbmUifX19XSxbMiwzLCJmIiwyLHsiY29sb3VyIjpbMCw2MCw2MF19LFswLDYwLDYwLDFdXSxbNSw0LCJmX1YiXSxbNCwzXSxbNSwyXSxbNSwzLCIiLDEseyJzdHlsZSI6eyJuYW1lIjoiY29ybmVyIn19XSxbOSw1LCIiLDIseyJsZXZlbCI6Miwic3R5bGUiOnsiaGVhZCI6eyJuYW1lIjoibm9uZSJ9fX1dLFs5LDgsImZeXFxzaGFycCIsMix7ImNvbG91ciI6WzAsNjAsNjBdfSxbMCw2MCw2MCwxXV0sWzExLDYsImZfe1xcYmV0YX0iXSxbMTAsNywiZl5cXHNoYXJwX1xcYWxwaGEiLDJdLFsxMCwxMSwiZl5cXGRvdWJsZXNoYXJwIiwwLHsiY29sb3VyIjpbMCw2MCw2MF0sInN0eWxlIjp7InRhaWwiOnsibmFtZSI6Im1vbm8ifX19LFswLDYwLDYwLDFdXSxbNyw4LCIiLDEseyJzdHlsZSI6eyJ0YWlsIjp7Im5hbWUiOiJtb25vIn19fV0sWzEwLDksIiIsMSx7InN0eWxlIjp7InRhaWwiOnsibmFtZSI6Im1vbm8ifX19XSxbMTEsNSwiIiwxLHsic3R5bGUiOnsidGFpbCI6eyJuYW1lIjoibW9ubyJ9fX1dLFs2LDQsIiIsMSx7InN0eWxlIjp7InRhaWwiOnsibmFtZSI6Im1vbm8ifX19XSxbOSwxXSxbOCwwXSxbMTEsNCwiIiwyLHsic3R5bGUiOnsibmFtZSI6ImNvcm5lciJ9fV0sWzEwLDgsIiIsMix7InN0eWxlIjp7Im5hbWUiOiJjb3JuZXIifX1dLFsxNCwxMywiUSJdLFsxMywxMiwiUCJdXQ==
		\begin{tikzcd}[ampersand replacement=\&]
			{\dual[\Simple{\cat C}]{\dual[\cat C]{\cat P}}} \&[-2ex]\&[-2ex] \alpha \& \textcolor{red}{(f^\sharp)^*\alpha} \& \textcolor{red}{f^*\beta} \& \beta \\
			{\Simple{\cat C}} \&\& \textcolor{red}{U \times X} \& \textcolor{red}{V \times X} \& {V \times X} \& {V \times Y} \\
			{\dual[\cat C]{\cat C}} \&\& X \& X \& \textcolor{red}{X} \& \textcolor{red}{Y}
			\arrow[Rightarrow, no head, from=3-3, to=3-4]
			\arrow[Rightarrow, no head, from=3-4, to=3-5]
			\arrow["f"', color={red}, from=3-5, to=3-6]
			\arrow["{f_V}", from=2-5, to=2-6]
			\arrow[from=2-6, to=3-6]
			\arrow[from=2-5, to=3-5]
			\arrow["\lrcorner"{anchor=center, pos=0.125}, draw=none, from=2-5, to=3-6]
			\arrow[Rightarrow, no head, from=2-4, to=2-5]
			\arrow["{f^\sharp}"', color={red}, from=2-4, to=2-3]
			\arrow["{f_{\beta}}", from=1-5, to=1-6]
			\arrow["{f^\sharp_\alpha}"', from=1-4, to=1-3]
			\arrow["{f^\doublesharp}", color={red}, tail, from=1-4, to=1-5]
			\arrow[tail, from=1-3, to=2-3]
			\arrow[tail, from=1-4, to=2-4]
			\arrow[tail, from=1-5, to=2-5]
			\arrow[tail, from=1-6, to=2-6]
			\arrow[from=2-4, to=3-4]
			\arrow[from=2-3, to=3-3]
			\arrow["\lrcorner"{anchor=center, pos=0.125}, draw=none, from=1-5, to=2-6]
			\arrow["\lrcorner"{anchor=center, pos=0.125, rotate=-90}, draw=none, from=1-4, to=2-3]
			\arrow["{\dual[\Simple{\cat C}]{\dual[\cat C]{P_2}}}", from=1-1, to=2-1]
			\arrow["{\dual[\cat C]{P_1}}", from=2-1, to=3-1]
		\end{tikzcd}
	\end{equation}

	Notice the tower of fibrations we used here is obtained quite naturally as a coproduct completion:
	\begin{equation}
		\begin{tikzcd}[ampersand replacement=\&,sep=5ex, column sep=3ex]
			{\SubFib{\cat C} \times_{\cat C} \Simple{\cat C}} \& {\SubFib{\cat C}} \\
			{\Simple{\cat C}} \& {\cat C} \\
			{\cat C}
			\arrow["\sub", from=1-2, to=2-2]
			\arrow["\times", from=2-1, to=2-2]
			\arrow["{\simple = P_1}"', from=2-1, to=3-1]
			\arrow[from=1-1, to=1-2]
			\arrow["{P_2}"', from=1-1, to=2-1]
			\arrow["\lrcorner"{anchor=center, pos=0.125}, draw=none, from=1-1, to=2-2]
		\end{tikzcd}
	\end{equation}

	Thus what we called `$\cat P$' above is really $\SubFib{\cat C} \times_{\cat C} \Simple{\cat C}$.

	\begin{definition}
		Let $\tower{P}$ be a tower of fibrations of length $n$.
		The morphisms in the total space of its iterated dual $\dual[n]{\tower{P}}$ are called \textbf{$\tower{P}$-dialenses}.
	\end{definition}

	Clearly dialenses of height $1$ are just $P$-lenses (Definition~\ref{def:lens}) and Dialectica (Definition~\ref{def:dc}) is a category of dialenses of height $2$.

	% \begin{remark}(\textbf{Why dialenses?})
	% 	The motivation for lenses given in Remark~\ref{rmk:why-lenses} generalizes naturally to dialenses of higher type. In fact, these can be seen to represent longer dialogues where at each stage, agents remember the previous ones.
	% \end{remark}

	\begin{example}
		As anticipated, optics can be constructed as dialenses of height $2$, if of a bizarre kind.
		First, recall from~\cite{braithwaite_fibre_2021} that every (strict, for simplicity) action $\bullet : \cat M \times \cat C \to \cat C$ of a monoidal category $\cat M$ on a category $\cat C$ can be repackaged as a 2-opfibration $\Para(\bullet) \to B\cat M$, where $\Para(\bullet)$ is the bicategory of $\cat M$-parametric $\cat C$-morphisms (see~\cite{Capucci_2022}) and the projection to $B\cat M$ is given by projection out the parameter of a morphism.
		Here let's ignore the 2-dimensional information and just consider this as an opfibration (because every 2-opfibration has an underlying 1-opfibration, cf.~\cite[Definition~2.1]{hermida_properties_1999}).

		Recall the 1-dimensional data of $\Para(\bullet)$: its objects are the same of $\cat C$ while its 1-cells are pairs $(m \in \cat M, f:m \bullet X \to X)$, composing in the obvious way.

		Opcartesian lifts for the functor $\Para(\bullet) \to B\cat M$ are constructed as follows:
		\begin{equation}
			% file:///home/jsb20179/data/software/quiver/src/index.html?q=WzAsNixbMCwwLCJcXFBhcmEoXFxidWxsZXQpIl0sWzAsMSwiQlxcY2F0IE0iXSxbMSwxLCJcXGFzdCJdLFsyLDEsIlxcYXN0Il0sWzEsMCwiWCJdLFsyLDAsIm0gXFxidWxsZXQgWCJdLFswLDEsIlAiLDJdLFsyLDMsIm0iLDJdLFs0LDUsIihtLCAxX3ttIFxcYnVsbGV0IFh9KSJdLFs1LDIsIikiLDMseyJsYWJlbF9wb3NpdGlvbiI6MCwic3R5bGUiOnsiYm9keSI6eyJuYW1lIjoibm9uZSJ9LCJoZWFkIjp7Im5hbWUiOiJub25lIn19fV1d
			\begin{tikzcd}[ampersand replacement=\&]
				{\Para(\bullet)} \& X \& {m \bullet X} \\
				{B\cat M} \& \ast \& \ast
				\arrow["P_\bullet"', from=1-1, to=2-1]
				\arrow["m"', from=2-2, to=2-3]
				\arrow["{(m, 1_{m \bullet X})}", from=1-2, to=1-3]
			\end{tikzcd}
		\end{equation}
		Further on, we'll abbreviate $(m, 1_{m \bullet X})$ with $m_X$.

		Then given two left $\cat M$-actegories $(\cat C,\bullet)$ and $(\cat D, \circ)$, one can construct the following tower of opfibrations:
		\begin{equation}
			% file:///home/jsb20179/data/software/quiver/src/index.html?q=WzAsNSxbMCwwLCJcXFBhcmEoXFxjaXJjKSBcXHRpbWVzX3tcXGNhdCBCTX0gXFxQYXJhKFxcYnVsbGV0KSJdLFswLDEsIlxcUGFyYShcXGJ1bGxldCkiXSxbMCwyLCJcXGNhdCBCTSJdLFsxLDEsIlxcY2F0IEJNIl0sWzEsMCwiXFxQYXJhKFxcY2lyYykiXSxbNCwzLCJcXGludCBcXGNpcmMiXSxbMSwzLCJcXGludCBcXGJ1bGxldCJdLFsxLDIsInBfMSA9IFxcaW50IFxcYnVsbGV0IiwyXSxbMCw0XSxbMCwxLCJwXzIiLDJdLFswLDMsIiIsMCx7InN0eWxlIjp7Im5hbWUiOiJjb3JuZXIifX1dXQ==
			\begin{tikzcd}[ampersand replacement=\&, sep=5ex]
				{\Para(\circ) \times_{B \cat M} \Para(\bullet)} \& {\Para(\circ)} \\
				{\Para(\bullet)} \& {B \cat M} \\
				{B\cat M}
				\arrow["{P_\circ}", from=1-2, to=2-2]
				\arrow["{P_\bullet}", from=2-1, to=2-2]
				\arrow["{P_\bullet}"', from=2-1, to=3-1]
				\arrow[from=1-1, to=1-2]
				\arrow["{P_\circ}"', from=1-1, to=2-1]
				\arrow["\lrcorner"{anchor=center, pos=0.125}, draw=none, from=1-1, to=2-2]
			\end{tikzcd}
		\end{equation}
		Then its iterated dual produces preoptics (Def.~\ref{def:preoptics}):
		\begin{equation}
			% file:///home/jsb20179/data/software/quiver/src/index.html?q=WzAsMTUsWzEsMiwiXFxhc3QiLFswLDYwLDYwLDFdXSxbMiwyLCJcXGFzdCIsWzAsNjAsNjAsMV1dLFszLDIsIlxcYXN0Il0sWzQsMiwiXFxhc3QiXSxbMSwxLCJYIl0sWzIsMSwibSBcXGJ1bGxldCBYIl0sWzMsMSwibSBcXGJ1bGxldCBYIixbMCw2MCw2MCwxXV0sWzQsMSwiWSIsWzAsNjAsNjAsMV1dLFs0LDAsIlYiXSxbMywwLCJWIixbMCw2MCw2MCwxXV0sWzEsMCwiVSJdLFsyLDAsIm0gXFxjaXJjIFUiLFswLDYwLDYwLDFdXSxbMCwyLCJCXFxjYXQgTSJdLFswLDEsIlxcZHVhbFtCXFxjYXQgTV17XFxQYXJhKFxcYnVsbGV0KX0iXSxbMCwwLCJcXE9wdGljX3tcXGNhdCBNfShcXGJ1bGxldCwgXFxjaXJjKSJdLFswLDEsIm0iLDAseyJjb2xvdXIiOlswLDYwLDYwXX0sWzAsNjAsNjAsMV1dLFsxLDIsIiIsMCx7ImxldmVsIjoyLCJzdHlsZSI6eyJoZWFkIjp7Im5hbWUiOiJub25lIn19fV0sWzIsMywiIiwwLHsibGV2ZWwiOjIsInN0eWxlIjp7ImhlYWQiOnsibmFtZSI6Im5vbmUifX19XSxbNCw1LCJtX1giXSxbNSw2LCIiLDAseyJsZXZlbCI6Miwic3R5bGUiOnsiaGVhZCI6eyJuYW1lIjoibm9uZSJ9fX1dLFs3LDYsImYiLDIseyJjb2xvdXIiOlswLDYwLDYwXX0sWzAsNjAsNjAsMV1dLFs5LDgsIjFfViIsMCx7ImxldmVsIjoyLCJzdHlsZSI6eyJoZWFkIjp7Im5hbWUiOiJub25lIn19fV0sWzEwLDExLCJtX1UiXSxbMTEsOSwiZl5cXHNoYXJwIiwwLHsiY29sb3VyIjpbMCw2MCw2MF19LFswLDYwLDYwLDFdXSxbMTMsMTIsIlxcZHVhbFtCXFxjYXQgTV17UF9cXGJ1bGxldH0iXSxbMTQsMTMsIlxcdmlldyJdXQ==
			\begin{tikzcd}[ampersand replacement=\&]
				{ \PreOptic_{\cat M}(\bullet, \circ)} \& U \& \textcolor{red}{m \circ U} \& \textcolor{red}{V} \& V \\
				{\dual[B\cat M]{\Para(\bullet)}} \& X \& {m \bullet X} \& \textcolor{red}{m \bullet X} \& \textcolor{red}{Y} \\
				{B\cat M} \& \textcolor{red}{\ast} \& \textcolor{red}{\ast} \& \ast \& \ast
				\arrow["m", color={red}, from=3-2, to=3-3]
				\arrow[Rightarrow, no head, from=3-3, to=3-4]
				\arrow[Rightarrow, no head, from=3-4, to=3-5]
				\arrow["{m_X}", from=2-2, to=2-3]
				\arrow[Rightarrow, no head, from=2-3, to=2-4]
				\arrow["f"', color={red}, from=2-5, to=2-4]
				\arrow["{1_V}", Rightarrow, no head, from=1-4, to=1-5]
				\arrow["{m_U}", from=1-2, to=1-3]
				\arrow["{f^\sharp}", color={red}, from=1-3, to=1-4]
				\arrow["{\dual[B\cat M]{P_\bullet}}", from=2-1, to=3-1]
				\arrow["\view", from=1-1, to=2-1]
			\end{tikzcd}
		\end{equation}
	\end{example}

	The idea of seeing optics as fibered by (or indexed by) their residual has appeared before in the literature, namely in \cite{braithwaite_fibre_2021,milewski_compound_2022,vertechi_dependent_2023,capucci_seeing_2022}.

	\section[Dialenses and Hofstra's Dial monad]{Dialenses and Hofstra's $\Dial$ monad}
	\label{sec:hofstra}
	In~\cite{hofstra_dialectica_2011}, Hofstra defined a monad on fibrations that builds Dialectica-like categories by simple sum-product completion:
	\begin{equation}
		\Dial(P) = \sSum(\sProd(P)) = \sSum(\dual{\sSum(\dual P)})
	\end{equation}
	where $P : \cat E \to \cat C$ is a fibration on $\cat C$ cartesian monoidal.
	Recall that simple sum completion $\sSum(P)$ is constructed like $\Fam(P)$ except $P$ gets pulled back along $\Simple(\cat C) \nto{\times} \cat C$ instead of $\dom$ (cf.~\eqref{eq:simple-dial-hofstra} below with~\eqref{eq:fam}).

	In $\Dial(P)$, objects are triples $(I, X, U : \cat C, \alpha : \cat P(I \times X \times U))$ and morphisms have \textbf{four} parts:
	\begin{eqalign}
		f_0 &\colon I \to J\\
		f &\colon I \times X \to Y\\
		f^\sharp &\colon I \times X \times V \to U\\
		\text{for } i : I, x:X,v:V,\  f^\doublesharp &\colon \alpha(i, x, f^\sharp(i,x,v)) \to \beta(f_0(i), f(i,x), v)
	\end{eqalign}
	If we ignore the duals for a moment (they can be put back later), we see this is actually given by a tower of \textbf{three} fibrations over $\cat C$:
	\begin{equation}
	\label{eq:simple-dial-hofstra}
		% file:///home/jsb20179/data/software/quiver/src/index.html?q=WzAsOSxbMCwyLCJcXGNhdCBDX1xccHJval5cXGRvd25hcnJvdyJdLFswLDMsIlxcY2F0IEMiXSxbMSwyLCJcXGNhdCBDIl0sWzEsMSwiXFxjYXQgQ19cXHByb2peXFxkb3duYXJyb3ciXSxbMiwxLCJcXGNhdCBDIl0sWzIsMCwiXFxjYXQgRSJdLFsxLDAsIlxcU3VtKHApIl0sWzAsMCwiXFxTdW0oXFxTdW0ocCkpIl0sWzAsMSwiXFxjYXQgQ19cXHByb2peXFxkb3duYXJyb3cgXFx0aW1lc197XFxjYXQgQ30gXFxjYXQgQ19cXHByb2peXFxkb3duYXJyb3ciXSxbNSw0LCJwIl0sWzYsM10sWzYsNV0sWzMsNCwiXFx0aW1lcyJdLFszLDIsIlxcY29kIiwxXSxbMCwyLCJcXHRpbWVzIl0sWzAsMSwicV8xID0gXFxjb2QiLDJdLFs3LDZdLFs3LDgsInFfMyIsMl0sWzgsMCwicV8yIiwyXSxbOCwzXSxbNiw0LCIiLDEseyJzdHlsZSI6eyJuYW1lIjoiY29ybmVyIn19XSxbNywxOSwiIiwxLHsibGV2ZWwiOjEsInN0eWxlIjp7Im5hbWUiOiJjb3JuZXIifX1dLFs4LDE0LCIiLDEseyJsZXZlbCI6MSwic3R5bGUiOnsibmFtZSI6ImNvcm5lciJ9fV1d
		\begin{tikzcd}[ampersand replacement=\&]
			{\sSum(\sSum(P))} \& {\sSum(P)} \& {\cat E} \\
			{\sSum(\simple)} \& {\Simple{\cat C}} \& {\cat C} \\
			{\Simple{\cat C}} \& {\cat C} \\
			{\cat C}
			\arrow["P", from=1-3, to=2-3]
			\arrow[from=1-2, to=2-2]
			\arrow[from=1-2, to=1-3]
			\arrow["\times", from=2-2, to=2-3]
			\arrow["\simple", from=2-2, to=3-2]
			\arrow[""{name=0, anchor=center, inner sep=0}, "\times", from=3-1, to=3-2]
			\arrow["{Q_1 = \simple}"', from=3-1, to=4-1]
			\arrow[from=1-1, to=1-2]
			\arrow["{Q_3}"', from=1-1, to=2-1]
			\arrow["{Q_2}"', from=2-1, to=3-1]
			\arrow[""{name=1, anchor=center, inner sep=0}, from=2-1, to=2-2]
			\arrow["\lrcorner"{anchor=center, pos=0.125}, draw=none, from=1-2, to=2-3]
			\arrow["\lrcorner"{anchor=center, pos=0.125}, draw=none, from=1-1, to=1]
			\arrow["\lrcorner"{anchor=center, pos=0.125}, draw=none, from=2-1, to=0]
		\end{tikzcd}
	\end{equation}

	One can see $\Dial(P)$ is obtained by dualizing the top two:
	\begin{equation}
		\Dial(P) = \dual[2]{(\sSum(\sSum(P)) \nlongto{Q_3} \Simple{\cat C} \times_{\cat C} \Simple{\cat C} \nlongto{Q_2} \Simple{\cat C})} \nlongto{Q_1} \cat C
	\end{equation}

	This readily carries over to the \emph{dependent Dialetica} construction
	\begin{equation*}
		\DepDial(P) := \Sum(\Prod(P))
	\end{equation*}
	of a fibration $P : \cat E \to \cat X$.\footnote{For a finitely complete category $\cat C$, we recover $\sDial(\cat C)$ as the fiber over the terminal object of $\DepDial(P)$, where $P : \cat E \to \cat X$ is taken to be $\Sum(!) : \Sum(\cat C) \to \Sum(1) \simeq \Set$, arising canonically from the free sum completion of $\cat C$, see~\cite{moss2022talk}.}
	Here, $\Sum(P)$ is the (unrestricted) \emph{sum completion} of $P$, which is the family construction of a fibration $P$, and $\Prod(P) := \dual{(\Sum(\dual{P}))}$ is the (unrestricted) \emph{product completion} of a fibration $P$. Thus, the iterated dual is applied to the following composite of fibrations:

	\begin{equation}
		\begin{tikzcd}
			{\Sum(\Sum(P))} & {\Sum(P)} & {\cat E} \\
			{\Sum(\cod)} & {\cat X^\downarrow} & {\cat X} \\
			{\cat X^\downarrow} & {\cat X} \\
			{\cat X}
			\arrow["P", from=1-3, to=2-3]
			\arrow[from=1-2, to=2-2]
			\arrow[from=1-2, to=1-3]
			\arrow["\dom", from=2-2, to=2-3]
			\arrow["\cod", from=2-2, to=3-2]
			\arrow[""{name=0, anchor=center, inner sep=0}, "\dom", from=3-1, to=3-2]
			\arrow["{{\cod}}"', from=3-1, to=4-1]
			\arrow[from=1-1, to=1-2]
			\arrow[from=1-1, to=2-1]
			\arrow[from=2-1, to=3-1]
			\arrow[""{name=1, anchor=center, inner sep=0}, from=2-1, to=2-2]
			\arrow["\lrcorner"{anchor=center, pos=0.125}, draw=none, from=1-2, to=2-3]
			\arrow["\lrcorner"{anchor=center, pos=0.125}, draw=none, from=1-1, to=1]
			\arrow["\lrcorner"{anchor=center, pos=0.125}, draw=none, from=2-1, to=0]
		\end{tikzcd}
	\end{equation}

	The elements of the total category of $\DepDial(P)$ are the dependent generalization of the simple case of $\Dial(P)$:
	\begin{eqalign}
		f_0 &\colon I \to J\\
		\text{for }i:I,\ f &\colon X(i) \to Y(f(i)) \\
		\text{for }i:I, x:X(i),\ f^\sharp &\colon V(f(i), f_0(i,x)) \to U(i,x)\\
		\text{for }i:I, x:X(i), y:Y(f(i), f_0(i,x)),\ f^\doublesharp &\colon \alpha(i,x,f_1(i,x,y)) \to \beta(f(i), f_0(i,x), y)
	\end{eqalign}

	\section*{Acknowledgement}
		We would like to thank Valeria de Paiva for her outstanding and sustained mentorship at the 2022 American Mathematical Society Mathematics Research Community (AMS~MRC) and beyond, as well as for providing feedback and corrections to a draft version of this text. Furthermore, we are grateful to numerous fruitful discussions with our other team members Charlotte Aten, Colin Bloomfield, Eric Bond, Joseph Dorta, Samantha Jarvis, Jérémie Koenig, Nelson Niu, and Jan Rooduijn. We also thank David Spivak for very insightful discussions of the material. We thank the three anonymous referees for valuable feedback and pointing out flaws in the preprint version of this text. This material is based upon work supported by the National Science Foundation under Grant Numbers
		DMS 1641020 and DMS 1916439. Any opinions, findings, and conclusions or recommendations expressed in this material
		are those of the authors and do not necessarily reflect the views of the National Science Foundation.

		JW is grateful for financial support by the US Army Research Office under MURI Grant W911NF-20-1-0082 hosted by the Department of Mathematics at Johns Hopkins University, and for hospitality and financial support provided by the Max Planck Institute for Mathematics, Bonn, Germany.
		
\bibliographystyle{entics}
\bibliography{bibliography}

\end{document}